\documentclass[12pt]{article}

\usepackage {amsmath, amssymb, amsthm}

\usepackage {tikz}
\usetikzlibrary{calc}

\usepackage{hyperref}
\hypersetup{
    colorlinks = true,
    citecolor = red,
}

\usepackage [margin=1.5in] {geometry}

\theoremstyle {definition}

\newtheorem {ex} {Example}

\theoremstyle {plain}
\newtheorem {lemma}{Lemma}
\newtheorem {thm}{Theorem}
\newtheorem {prop}[thm]{Proposition}

\newcommand{\C}{{\mathbb C}}

\newcommand{\G}{{\mathcal G}}
\newcommand{\be}{\begin{equation}}
\newcommand{\ee}{\end{equation}}

\newcommand{\R}{{\mathbb R}}
\newcommand{\Z}{{\mathbb Z}}

\newcommand{\GL}{\mathrm{GL}}

\newcommand{\XX}{\mathcal{X}}
\newcommand{\YY}{\mathcal{Y}}

\newcommand{\OM}{\Omega}
\newcommand{\ODS}{\Omega^{\text{DS}}}

\newcommand{\old}[1]{}

\begin {document}
\title {Eigenvalues of matrix products}
\author{Richard Kenyon\footnote{Department of Mathematics, Yale University, New Haven; richard.kenyon at yale.edu} \and
Nicholas Ovenhouse\footnote{Department of Mathematics, Yale University, New Haven; nicholas.ovenhouse at yale.edu}}
\date{}
\maketitle

\begin{abstract}
We study pairs of matrices $A,B\in \GL_n(\C)$
such that the eigenvalues of $A$, of $B$ and of the product $AB$ are specified in advance. 
We show that the space of such pairs $(A,B)$ under simultaneous conjugation  
has dimension $(n-1)(n-2)$, and give an explicit parameterization.

More generally let $\Sigma$ be a surface of genus $g$ with $k$ punctures.
We find a parameterization of the space $\Omega_{g,k,n}$ of 
flat $\GL_n(\C)$-structures on $\Sigma$ whose holonomies
around the punctures have prescribed
eigenvalues.

We show furthermore that, for $3\le k\le 2g+6$ (or $3\le k\le 9$ if $g=1$, or $3\le k$ if $g=0$), the space $\Omega_{g,k,n}$ has an explicit symplectic structure and an associated Liouville integrable system,
equivalent to a leaf of a Goncharov-Kenyon dimer integrable system. 
\end{abstract}

\section{Introduction}

The \emph{Deligne-Simpson problem} (DSP), see \cite{Kostov}, is to find matrices $M_1,\dots,M_k$ 
whose product is the identity 
and which have prescribed conjugacy classes in $\GL_n(\C)$.
A more general setting is as follows.
Given an oriented surface $\Sigma$ of genus $g$ with $k$ boundary components,
describe the space 
$\ODS_{g,k,n}(C)$ of representations $\pi_1(\Sigma)\to\GL_n(\C)$
with prescribed conjugacy classes $\{C_i\}_{i=1,\dots,k}$ 
for the peripheral curves around each boundary component.
This problem arises from the study of systems of singular linear ODEs on surfaces. 

The simplest case is when $\Sigma$ is a $3$-holed sphere (the $k=3$ case of the DSP).
In this case,
this problem can be restated as finding matrices $A,B\in\GL_n(\C)$, for which the conjugacy classes of $A,B$ and $AB$ are prescribed. 

We deal here with the simpler version of prescribing the eigenvalues only.
Let $\OM_{g,k,n}(\lambda)$ be the space of representations $\pi_1(\Sigma)\to\GL_n(\C)$, up to conjugacy, where $
\lambda=(\lambda_1,\dots,\lambda_k)$ are the sets of eigenvalues for the monodromies of the $k$ peripheral curves.
That is, $\lambda_i := \{\lambda_{i,1}, \lambda_{i,2}, \dots, \lambda_{i,n}\}$ are the eigenvalues of the monodromy around boundary component $i$.
Note that if each $\lambda_i$ consists of distinct values, this is just the DSP,
that is,
$\OM_{g,k,n}(\lambda)=\ODS_{g,k,n}(\lambda)$ when each $\lambda_i$ consists of distinct eigenvalues. When some individual eigenvalues are not distinct,
$\OM_{g,k,n}(\lambda)$ is a union of subspaces of various dimensions consisting of $\ODS_{g,k,n}(C)$ for the various conjugacy classes $C$ having the same eigenvalues as $\lambda$. 

Our main result is a Laurent parameterization of the space $\OM_{g,k,n}(\lambda)$,  
$$ \Psi:(\C^*)^N\to \OM_{g,k,n}(\lambda), $$
where $N=2(g-1)n^2+(n^2-n)k+2$ is the dimension of $\OM_{g,k,n}(\lambda)$. That is, we give explicit matrices representing a connection in $\OM_{g,k,n}(\lambda)$
whose entries are Laurent polynomials in $N$ parameters. When $(g,k)=(0,3)$ and $n$ is odd, and $\alpha,\beta,\gamma$ consist in positive reals, 
there is a positive Laurent parameterization (with coefficients $1$) $\Psi:(\R_+)^N\to\OM_{0,3,n}^+(\lambda)$ onto a certain ``positive" subvariety
$\OM_{0,3,n}^+(\lambda)\subset\OM_{0,3,n}(\lambda)$, see Section \ref{positive}.

\begin{figure}[htbp]
\begin{center}\includegraphics[width=2in]{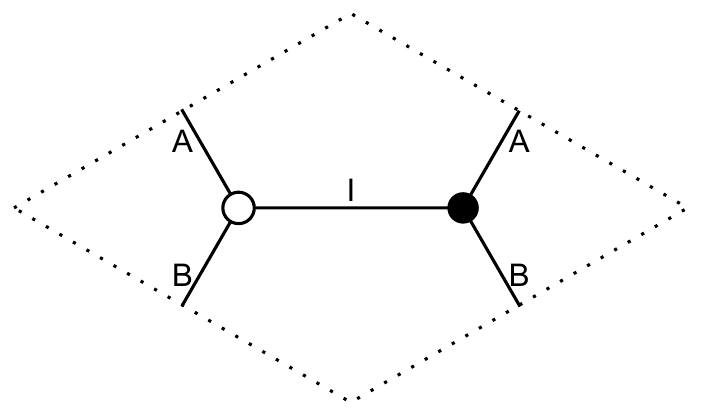}\end{center}
\caption{\label{G0}The graph $\G$ on the torus (the torus is obtained by identifying opposite sides of the rhombus) has two vertices and three edges.
The parallel transports along the three edges from black to white are $A,B,I$ as shown.}
\end{figure}

For $(g,k)=(0,3)$ a similar problem, Horn's problem, asks to prescribe the \emph{singular values} of three matrices $A,B$ and $AB$; 
this was solved by Knutson and Tao \cite{KT}, based on earlier work of Klyachko \cite{Kly},
by finding a surprising connection with ``hives" or tropical degree-$n$ curves.
Our method uses a similar (but pre-tropical) setup. 
It arises from studying $\GL_n$-local systems on certain graphs on surfaces
(see Figure \ref{G0} for the $(g,k)=(0,3)$ case), and is inspired by the work of Fock and Goncharov \cite{FG}. 

The basic idea behind the construction is that there is an isomorphism between $\GL_n$
structures on a bipartite surface graph $\G$, and $\C^*$-structures on a related graph $\G_n$. 
This idea 
originated in \cite{FG} who coordinatized the space of (decorated) $\GL_n$ structures
with explicit scalar quantities on the related graphs $\G_n$. See also \cite{KO24}.

If $3\le k\le 2g+6$, or $g=1$ and $3\le k\le 9$, or $g=0$ and $k\ge 3$,
we show in Section \ref{intscn} that 
the space $\OM_{g,k,n}(\lambda)$ can be identified
with a symplectic leaf of a dimer integrable system as defined in \cite{GK}. That is, there is a certain nondegenerate Poisson structure on $\OM_{g,k,n}(\lambda)$ and a set of 
explicit commuting Hamiltonian
functions on it forming a Liouville integrable system. 
When $k$ is outside the above bounds we don't know if there is such a correspondence.

\subsection{Organization}
In Sections \ref{conjugationscn}, \ref{pathscn}, and \ref{paramscn} we deal with the simplest nontrivial case of $(g,k)=(0,3)$. 
This is the case of parameterizing pairs $(A,B)$ of matrices in $\GL_n$ with prescribed eigenvalues
for $A,B$ and $AB$. We discuss the general case of punctured surfaces in Section \ref{genscn}.
The integrability results are discussed in Section \ref{intscn}.

\bigskip

\noindent{\bf Acknowledgments.} 
This research was supported by NSF grant DMS-1940932 and the Simons Foundation grant 327929. 
We thank Daniel Douglas and Sergey Fomin for discussions.

\section{Conjugation}\label{conjugationscn}

We replace $A$ with $A^{-1}$ and parameterize pairs $(A,B)$ with prescribed eigenvalues for $A,B$ and $BA^{-1}$. See Theorem \ref{30casethm}.

Given three multisets (sets allowing repeated elements)  $\alpha,\beta,\gamma\subset\C^*$, each of size $n$:
$$\alpha=\{\alpha_1,\dots,\alpha_n\},~~\beta=\{\beta_1,\dots,\beta_n\},~~\gamma=\{\gamma_1,\dots,\gamma_n\},$$
let $\Omega=\OM_{0,3,n}(\alpha,\beta,\gamma)$ be the space of pairs of matrices $A,B\in\GL_n(\C)$, up to simultaneous conjugation,
with eigenvalues $\alpha$ and $\beta$ respectively, 
and with product $BA^{-1}$ having eigenvalues $\gamma$.
Note that simultaneous conjugation of $A,B$ does not change any of $\alpha, \beta, \gamma$.

Let $f_1\subseteq \dots \subseteq f_n$ be a complete flag invariant under $A$: this means for each $i$,
$f_i$ is subspace of $\C^n$ of dimension $i$ and $Af_i\subseteq f_i$. Likewise let 
$g_1\subseteq \dots \subseteq g_n$ be an invariant flag for $B$.
Further suppose that the two flags are \emph{transverse}; that is, $\dim(f_i \cap g_j) = \max(0, i+j-n)$; this will be true for generic $A,B$.

\begin {lemma}
    Let $M$ be a matrix whose $k$th column spans the 1-dimensional space $f_{n+1-k}\cap g_{k}$.
    Then $M^{-1}AM$ is lower triangular and $M^{-1}BM$ is upper triangular.
\end {lemma}
\begin {proof}
    Let $M_k$ be the $k$th column of the matrix $M$ (which spans $f_{n+1-k} \cap g_k$). 
    Since $M_k \in g_k$, and the $M_k$'s are linearly independent, the flag defined by the columns (in left-to-right order)
    \[ \left< M_1 \right> \subset \left< M_1, M_2 \right> \subset \left< M_1, M_2, M_3 \right> \subset \cdots \]
    is the same as the flag $g$. Since $M^{-1}BM$ is the coordinate expression of $B$ in the $M_k$ basis, and since $g$ is invariant under $B$,
    we see that $M^{-1}BM$ is upper-triangular. Similarly, the flag $f$ (invariant under $A$) is the same as the flag generated by the right-to-left ordering of the columns of $M$:
    \[ \left< M_n \right> \subset \left< M_n, M_{n-1} \right> \subset \left< M_n, M_{n-1}, M_{n-2} \right> \subset \cdots, \]
    and therefore $M^{-1}AM$ is lower-triangular.
\end {proof}

Since for generic $A,B$ the flags
$f,g$ are transverse, in the nongeneric case we can, by taking limits, 
still assume $A,B$ are lower and upper triangular.
This shows that without loss of generality we may assume that 
$A,B$ are lower and upper triangular, respectively. 

We may further conjugate by diagonal matrices (without loss of generality of determinant $1$), since these preserve upper and lower triangularity. 
The resulting space of matrix pairs with eigenvalues $\alpha,\beta$
up to conjugacy is then of dimension $n(n-1)-(n-1)=(n-1)^2$. Imposing $n-1$ eigenvalue conditions on the product (one condition is
determined, since $\det BA^{-1} = \det B/\det A$)
reduces the dimension of the space of pairs in $\Omega$ to $(n-1)(n-2)$.

Note that by choosing different flags for $A$ and/or for $B$ we can reorder the diagonal elements of $M^{-1}AM$ and of 
$M^{-1}BM$ independently. Our parameterization will depend on this choice, which is, equivalently, a choice of order for the eigenvalues of $A$ and of $B$.
It will also depend on a choice of order for $\gamma$, the eigenvalues of $BA^{-1}$.

\section {Paths in the honeycomb}\label{pathscn}

\begin {figure}[h]
\centering
\begin {tikzpicture}[scale=0.8, every node/.style={scale=0.8}]
    \foreach \x in {0,1,2,3} {
        \foreach \y in {0,1,2,3} {
        \draw ($({1.5*(\y+\x)},{0.5*sqrt(3)*(\y-\x)})$) -- ($({1.5*(\y+\x)},{0.5*sqrt(3)*(\y-\x)}) + (60:1)$);
        \draw ($({1.5*(\y+\x)},{0.5*sqrt(3)*(\y-\x)})$) -- ($({1.5*(\y+\x)},{0.5*sqrt(3)*(\y-\x)}) + (-60:1)$);
        \draw ($({1.5*(\y+\x)},{0.5*sqrt(3)*(\y-\x)})$) -- ($({1.5*(\y+\x)},{0.5*sqrt(3)*(\y-\x)}) + (180:1)$);
        \draw ($({1.5*(\y+\x)},{0.5*sqrt(3)*(\y-\x)}) + (-1,0)$) -- ($({1.5*(\y+\x)},{0.5*sqrt(3)*(\y-\x)}) + (-1,0) + (-120:1)$);
        \draw ($({1.5*(\y+\x)},{0.5*sqrt(3)*(\y-\x)}) + (-1,0)$) -- ($({1.5*(\y+\x)},{0.5*sqrt(3)*(\y-\x)}) + (-1,0) + (120:1)$);

        \draw ($({1.5*(\y+\x)},{0.5*sqrt(3)*(\y-\x)}) + (30:0.6)$) node {\scriptsize $a_{\x\y}$};
        \draw ($({1.5*(\y+\x)},{0.5*sqrt(3)*(\y-\x)}) + (-30:0.6)$) node {\scriptsize $b_{\x\y}$};
        }
    }
    \foreach \x in {0,1,2,3} {
        \foreach \y in {0,1,2,3} {
        \draw[fill=black] ($({1.5*(\y+\x)},{0.5*sqrt(3)*(\y-\x)})$) circle (0.08);
        \draw[fill=white] ($({1.5*(\y+\x)},{0.5*sqrt(3)*(\y-\x)}) + (-1,0)$) circle (0.08);
        }
    }
\end {tikzpicture}
\caption{The honeycomb graph $H$.}
\label{honeycomb}
\end {figure}
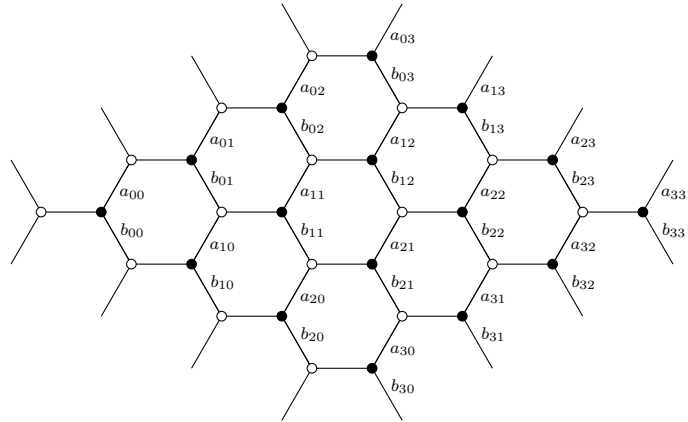

Let $H$ be the honeycomb graph (Figure \ref{honeycomb}). It is bipartite; orient its 
edges from black to white. 
We put edge weights $1$ on horizontal edges, $a_v$ on NE and $b_v$ on SE edges out
of black vertex $v$. These weights define a $\C^*$-connection $\phi$ on $H$, as follows. 
Associate a $1$-dimensional $\C$-vector space $\C_v\cong\C$ to each vertex $v$. 
Then define the isomorphism $\phi_{bw}:\C_b\to\C_w$ as multiplication by the edge weight $c_{bw}$. The isomorphism associated to the reversed edge is the inverse: $\phi_{wb} = \phi_{bw}^{-1}$. 

    For any vertex $v$, and non-zero scalar $t$, the associated \emph{gauge transformation} is the operation of multiplying
    the weights of all edges incident to $v$ by $t$. The group $(\Bbb{C}^*)^{V}$ generated by all such transformations is called the \emph{gauge group}.
    Graph connections are typically considered up to gauge equivalence, since equivalent connections have the same monodromies around loops. Our convention of setting all horizontal
    edge weights equal to $1$ is therefore no loss in generality; we may always obtain such weights with a particular choice of gauge.

\subsection{Monotone paths}

We consider three families of paths on $H$, called east (E) paths, northwest (NW) paths, and southwest (SW) paths.
East paths are paths which follow the orientation of $H$ once we reverse orientation on all horizontal edges. They are thus monotone left-to-right paths.
SW paths are paths which follow the orientation of $H$ once we reverse orientation on all NE edges.
NW paths are paths which follow the orientation of $H$ once we reverse orientation on all SE edges. 

\subsection{Paths in triangles}

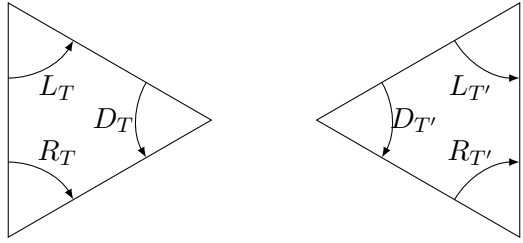
\begin{figure}[h]
\centering
\begin {tikzpicture}
    \def\r{1.8}
    \draw (0:\r) -- (120:\r) -- (-120:\r) -- cycle;
    \draw[-latex] ($(120:\r) - (0,1)$) arc (-90:-30:1); 
    \draw[-latex] ($(-120:\r) + (0,1)$) arc (90:30:1); 
    \draw[-latex] ($(0:\r) + (150:1)$) arc (150:210:1); 
    \draw (120:0.5)  node {\small $L_T$};
    \draw (-120:0.5) node {\small $R_T$};
    \draw (0:0.5)    node {\small $D_T$};

    \begin {scope}[shift={(5,0)}]
        \draw (180:\r) -- (-60:\r) -- (60:\r) -- cycle;
        \draw[-latex] ($(60:\r) + (-150:1)$) arc (-150:-90:1); 
        \draw[-latex] ($(-60:\r) + (150:1)$) arc (150:90:1); 
        \draw[-latex] ($(180:\r) + (30:1)$)  arc (30:-30:1); 
        \draw (60:0.5)  node {\small $L_{T'}$};
        \draw (-60:0.5) node {\small $R_{T'}$};
        \draw (180:0.5) node {\small $D_{T'}$};
    \end {scope}
\end {tikzpicture}
\caption{The triangles $T$ (left) and $T'$ (right), and the orientations of their associated transport matrices.}
\label{TT}
\end{figure}

We consider two types of triangular regions $T$ and $T'$, 
oriented as in  Figure \ref{TT}.

\begin {figure}[h]
\centering
\begin {tikzpicture}[scale=0.8]

    \draw (0:4) -- (120:4) -- (-120:4) -- cycle;

    \foreach \x/\y[evaluate={\X=int(\x+3); \Y=int(2-\y);}] in {0/0, 1/0, -1/0, -2/1, -1/1, 0/1, 1/1, 0/-1, 1/-1, 1/-2} {
        \draw ($({1.5*(\y+\x)},{0.5*sqrt(3)*(\y-\x)})$) -- ($({1.5*(\y+\x)},{0.5*sqrt(3)*(\y-\x)}) + (60:1)$);
        \draw ($({1.5*(\y+\x)},{0.5*sqrt(3)*(\y-\x)})$) -- ($({1.5*(\y+\x)},{0.5*sqrt(3)*(\y-\x)}) + (-60:1)$);
        \draw ($({1.5*(\y+\x)},{0.5*sqrt(3)*(\y-\x)})$) -- ($({1.5*(\y+\x)},{0.5*sqrt(3)*(\y-\x)}) + (180:1)$);

        \draw ($({1.5*(\y+\x)},{0.5*sqrt(3)*(\y-\x)}) + (30:0.6)$) node {\scriptsize $a_{\X\Y}$};
        \draw ($({1.5*(\y+\x)},{0.5*sqrt(3)*(\y-\x)}) + (-30:0.6)$) node {\scriptsize $b_{\X\Y}$};
    }
    \foreach \x/\y in {0/0, 1/0, -1/0, -2/1, -1/1, 0/1, 1/1, 0/-1, 1/-1, 1/-2} {
        \draw[fill=black] ($({1.5*(\y+\x)},{0.5*sqrt(3)*(\y-\x)})$) circle (0.08);
    }
    \foreach \x/\y in {0/0, 1/0, -1/0, 0/-1, 1/-1, 1/-2} {
        \draw[fill=white] ($({1.5*(\y+\x)+0.5},{0.5*sqrt(3)*(\y-\x+1)})$) circle (0.08);
    }

    \draw ($(-120:4) + (-1,{sqrt(3)/2})$) node {$4$};
    \draw ($(-120:4) + (-1,{3*sqrt(3)/2})$) node {$3$};
    \draw ($(-120:4) + (-1,{5*sqrt(3)/2})$) node {$2$};
    \draw ($(-120:4) + (-1,{7*sqrt(3)/2})$) node {$1$};

    \draw ($(120:4) + {sqrt(3)/2}*(-30:1) + (60:1)$) node {$1$};
    \draw ($(120:4) + {3*sqrt(3)/2}*(-30:1) + (60:1)$) node {$2$};
    \draw ($(120:4) + {5*sqrt(3)/2}*(-30:1) + (60:1)$) node {$3$};
    \draw ($(120:4) + {7*sqrt(3)/2}*(-30:1) + (60:1)$) node {$4$};

    \draw ($(0:4) + {sqrt(3)/2}*(-150:1) + (-60:1)$) node {$1$};
    \draw ($(0:4) + {3*sqrt(3)/2}*(-150:1) + (-60:1)$) node {$2$};
    \draw ($(0:4) + {5*sqrt(3)/2}*(-150:1) + (-60:1)$) node {$3$};
    \draw ($(0:4) + {7*sqrt(3)/2}*(-150:1) + (-60:1)$) node {$4$};

    \begin {scope}[shift={(9,0)}]
        \draw (60:4) -- (180:4) -- (-60:4) -- cycle;

        \foreach \x/\y in {0/0, 1/0, 2/0, 0/1, 1/1, 0/2} {
            \draw ($({1.5*(\y+\x)-2},{0.5*sqrt(3)*(\y-\x)})$) -- ($({1.5*(\y+\x)-2},{0.5*sqrt(3)*(\y-\x)}) + (60:1)$);
            \draw ($({1.5*(\y+\x)-2},{0.5*sqrt(3)*(\y-\x)})$) -- ($({1.5*(\y+\x)-2},{0.5*sqrt(3)*(\y-\x)}) + (-60:1)$);
            \draw ($({1.5*(\y+\x)-2},{0.5*sqrt(3)*(\y-\x)})$) -- ($({1.5*(\y+\x)-2},{0.5*sqrt(3)*(\y-\x)}) + (180:1)$);

            \draw ($({1.5*(\y+\x)-2},{0.5*sqrt(3)*(\y-\x)}) + (30:0.6)$) node {\scriptsize $a_{\x\y}$};
            \draw ($({1.5*(\y+\x)-2},{0.5*sqrt(3)*(\y-\x)}) + (-30:0.6)$) node {\scriptsize $b_{\x\y}$};
        }

        \foreach \x/\y in {3/0, 2/1, 1/2, 0/3} {
            \draw ($({1.5*(\y+\x)-3},{0.5*sqrt(3)*(\y-\x)})$) --++ (1,0);
        }
        \foreach \x/\y in {0/0, 1/0, 2/0, 3/0} {
            \draw ($({1.5*(\y+\x)-3},{0.5*sqrt(3)*(\y-\x)})$) --++ (-120:1);
        }
        \foreach \x/\y in {0/0, 0/1, 0/2, 0/3} {
            \draw ($({1.5*(\y+\x)-3},{0.5*sqrt(3)*(\y-\x)})$) --++ (120:1);
        }

        \foreach \x/\y in {0/0, 1/0, 2/0, 0/1, 1/1, 0/2} {
            \draw[fill=black] ($({1.5*(\y+\x)-2},{0.5*sqrt(3)*(\y-\x)})$) circle (0.08);
        }
        \foreach \x/\y in {0/0, 1/0, 2/0, 3/0, 0/1, 1/1, 2/1, 0/2, 1/2, 0/3} {
            \draw[fill=white] ($({1.5*(\y+\x)-3},{0.5*sqrt(3)*(\y-\x)})$) circle (0.08);
        }

        \draw ($(60:4) + {sqrt(3)/2}*(-90:1)   + (1,0)$) node {$1$};
        \draw ($(60:4) + {3*sqrt(3)/2}*(-90:1) + (1,0)$) node {$2$};
        \draw ($(60:4) + {5*sqrt(3)/2}*(-90:1) + (1,0)$) node {$3$};
        \draw ($(60:4) + {7*sqrt(3)/2}*(-90:1) + (1,0)$) node {$4$};

        \draw ($(180:4) + {sqrt(3)/2}*(30:1)   + (120:1)$) node {$4$};
        \draw ($(180:4) + {3*sqrt(3)/2}*(30:1) + (120:1)$) node {$3$};
        \draw ($(180:4) + {5*sqrt(3)/2}*(30:1) + (120:1)$) node {$2$};
        \draw ($(180:4) + {7*sqrt(3)/2}*(30:1) + (120:1)$) node {$1$};

        \draw ($(180:4) + {sqrt(3)/2}*(-30:1)   + (-120:1)$) node {$1$};
        \draw ($(180:4) + {3*sqrt(3)/2}*(-30:1) + (-120:1)$) node {$2$};
        \draw ($(180:4) + {5*sqrt(3)/2}*(-30:1) + (-120:1)$) node {$3$};
        \draw ($(180:4) + {7*sqrt(3)/2}*(-30:1) + (-120:1)$) node {$4$};
    \end {scope}
\end {tikzpicture}
\caption{The honeycomb graphs $H_T$ (left) and $H_{T'}$ (right), pictured for $n=4$.}
\label{Tn}
\end {figure}
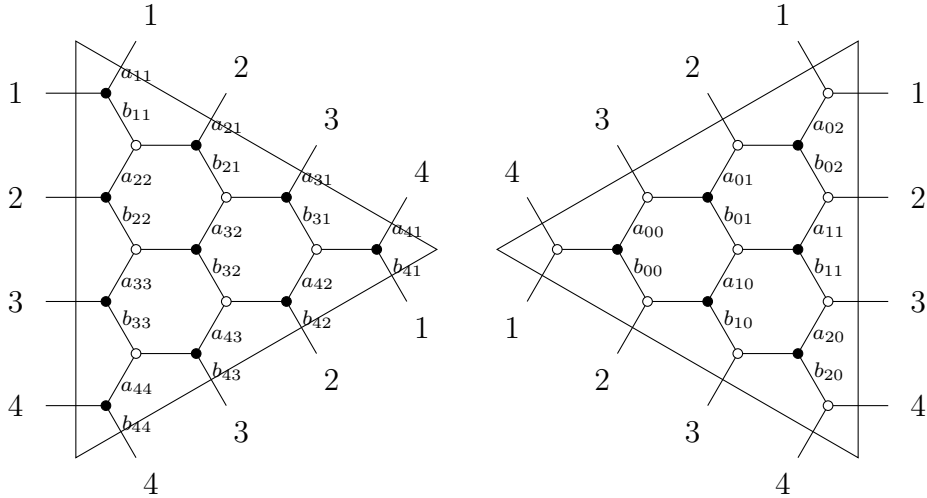

In $T$, for a given $n$ we define a honeycomb graph $H_T$ as shown in Figure \ref{Tn}, left.
It is the subgraph of the honeycomb $H$ which fits inside the triangle $T$ with side length $n$.
Edges of $H_T$ crossing the edges of each triangle are indexed from top to bottom as shown.
The graph $H_T$ has edge weights inherited from $H$, with $1$ on horizontal edges, and weights $a_v,b_v$ on the NE and
SE edges out of a black vertex $v$. 
In the figure, we have labeled the vertices by pairs $ij$, depending on their coordinates in the plane.

Associated to $H_T$ are two $n\times n$ matrices $L=L_T$ and $R=R_T$, 
where $L,R$ stand for left-turn and right-turn, respectively. The matrix $L$ is lower triangular and $R$ is upper triangular;
they are defined as follows. The entry $L_{ij}$ is the weighted sum of east paths from the $j$th edge
on the left side of $H_T$ to the $i$th edge on the top-right side of $H_T$. The weight of a path
is the holonomy of the $\C^*$ connection along the path: in this case it is just the product of the edge weights along the path, 
including the last edge. Similarly $R_{ij}$ is the weighted sum of holonomies of
east paths from the $j$th edge
on the left side of $H_T$ to the $i$th edge on the bottom-right side of $H_T$, including the last edge.

\begin {ex}
    The left-turn matrix for $n=3$ (corresponding to Figure \ref{Tn} with the bottom-right row of vertices removed) is given by
    \[ 
        L_T = \begin{pmatrix} 
            a_{11}              & 0                                        & 0                   \\
            b_{11}a_{21}        & a_{22}a_{21}                             & 0                   \\
            b_{11}b_{21}a_{31}  & a_{22}b_{21}a_{31} + b_{22}a_{32}a_{31}  & a_{33}a_{32}a_{31}  
        \end{pmatrix} 
    \]
\end {ex}

Likewise for the triangle $T'$ there is a honeycomb graph $H_{T'}$ (Figure \ref{Tn}, right). We define a 
corresponding (upper triangular) left-turn matrix $L'=L_{T'}$ and (lower triangular) right-turn matrix $R'=R_{T'}$. The entry $L'_{ij}$ is the weighted sum of east paths from the $j$th edge on the upper left side of $H_{T'}$ to the $i$th
edge on the right side of $H_{T'}$ (the weight does \emph{not} include the weight of the first edge).
Likewise $R'_{ij}$ is the weighted sum of east paths from the $j$th edge on the lower left side of $H_{T'}$ to the $i$th
edge on the right side of $H_{T'}$, with weight not including the first edge's weight.

\begin {ex}
    The right turn matrix for $n=4$ in $T'$ (see Figure \ref{Tn}, right) is 
    \[
        R_{T'} = \begin{pmatrix}
            a_{00}a_{01}a_{02}                                           & 0                           & 0      & 0 \\
            a_{00}a_{01}b_{02} + a_{00}b_{01}a_{11} + b_{00}a_{10}a_{11} & a_{10}a_{11}                & 0      & 0 \\
            a_{00}b_{01}b_{11} + b_{00}a_{10}b_{11} + b_{00}b_{10}a_{20} & b_{10}a_{20} + a_{10}b_{11} & a_{20} & 0 \\
            b_{00}b_{10}b_{20}                                           & b_{10}b_{20}                & b_{20} & 1 
        \end{pmatrix}
    \]
\end {ex}

A third matrix $D=D_T=R_TL_T^{-1}$ is associated to $T$,
and a third matrix $D'=D_{T'}=(R_{T'})^{-1}L_{T'}$ is associated to $T'$. 
The following is a combinatorial interpretation for the entries of $D$ and $D'$, which
appeared in a slightly different (but equivalent) form in a paper by Chekhov and Shapiro.

\begin{lemma} [\cite{CS}, Theorem 2.14] \label{D-paths}
    The entry $D_{ij}$
    is $(-1)^{j+n}$ times the weighted sum of SW paths from the $j$th edge on the top of $T$ to the $i$th edge on the bottom of $T$, including start and finish edges. 
    The entry $D'_{ij}$
    is $(-1)^{i+1}$ times the weighted sum of  SW paths from the $j$th edge on the top of $T'$ to the $i$th
    edge on the bottom of $T'$, not including start or finish edges. 
    The matrix $D$ is lower antidiagonal, that is $D_{i,j}=0$ for $i+j< n+1$, and $D'$ is upper antidiagonal: $D'_{i,j}=0$ for $i+j> n+1$.
\end{lemma}

\begin {ex}
    For $n=3$, we have
    \[ 
        D = \begin{pmatrix} 
            0 & 0 & a_{31}^{-1}b_{31} \\
            0 & -a_{21}^{-1}b_{21}a_{32}^{-1}b_{32} & a_{31}^{-1}a_{32}^{-1}b_{32} \\
            a_{11}^{-1}b_{11}a_{22}^{-1}b_{22}a_{33}^{-1}b_{33} & -a_{21}^{-1}a_{22}^{-1}b_{22}a_{33}^{-1}b_{33} - a_{21}^{-1}b_{21}a_{32}^{-1}a_{33}^{-1}b_{33} & a_{31}^{-1}a_{32}^{-1}a_{33}^{-1}b_{33} 
        \end{pmatrix}
    \]
\end {ex}

\subsection{The $A$ and $B$ matrices}
Given a positive integer $n$, consider the quadrilateral formed by adjacent triangles $T$ and $T'$ of size $n$.
We assign $a$ and $b$ edge weights as above inside the two triangles, and then extend this periodically
to the whole infinite honeycomb graph $H$. Equivalently, we identify opposite sides of this quadrilateral, and think of the union of $H_T$ and $H_{T'}$
as a graph $\G=\G_n$ drawn on a torus.

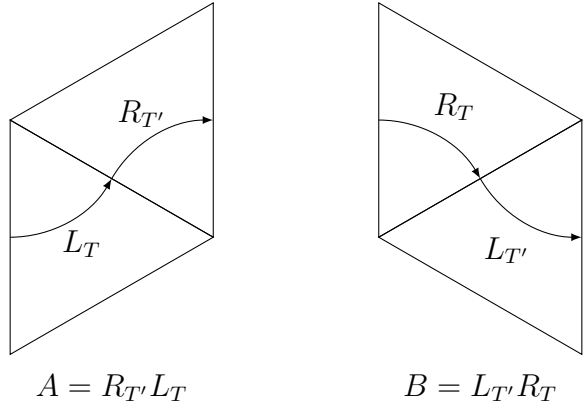
\begin{figure}[h]
\centering
\begin {tikzpicture}
    \def\r{1.8}
    \draw (0:\r) -- (120:\r) -- (-120:\r) -- cycle;
    \draw[-latex] ($(120:\r) - sqrt(3)*0.5*\r*(0,1)$) arc (-90:-30:{sqrt(3)*0.5*\r}); 
    \draw (-60:0.1) node {$L_T$};

    \begin {scope}[shift={(60:\r)}]
        \draw (180:\r) -- (-60:\r) -- (60:\r) -- cycle;
        \draw[-latex] ($(-60:\r) + sqrt(3)*0.5*\r*(150:1)$) arc (150:90:{sqrt(3)*0.5*\r}); 
        \draw (120:0.1) node {$R_{T'}$};
    \end {scope}

    \draw (0.25*\r, -2) node {$A = R_{T'}L_T$};

    \begin {scope}[shift={($(4,0) + (60:\r)$)}]
    \draw (0:\r) -- (120:\r) -- (-120:\r) -- cycle;
    \draw[-latex] ($(-120:\r) + sqrt(3)*0.5*\r*(0,1)$) arc (90:30:{sqrt(3)*0.5*\r}); 
    \draw (60:0.2) node {$R_T$};

    \begin {scope}[shift={(-60:\r)}]
        \draw (180:\r) -- (-60:\r) -- (60:\r) -- cycle;
        \draw[-latex] ($(60:\r) + sqrt(3)*0.5*\r*(-150:1)$) arc (-150:-90:{sqrt(3)*0.5*\r}); 
        \draw (-120:0.2) node {$L_{T'}$};

        \draw (-0.25*\r, -2) node {$B = L_{T'}R_T$};
    \end {scope}
    \end {scope}
\end {tikzpicture}
\caption{The networks associated to the $A$ and $B$ matrices.}
\label{BA_matrices}
\end{figure}

Given this initial data, we will define two matrices $A$ and $B$, one upper-triangular and one lower-triangular. 
They are given as the compositions of the $L$ and $R$ matrices defined above:
\[ A := R_{T'} L_T, \quad \quad B := L_{T'} R_T \]
By construction, $A$ is lower triangular (since both $R_{T'}$ and $L_T$ are lower-triangular), and similarly $B$ is upper triangular.
From the combinatorial description of the $R$ and $L$ matrices in terms of paths, it is easy to see that the entries of both
$A$ and $B$ are generating functions for east paths in networks formed by concatenating $H_T$ and $H_{T'}$ (using slightly
different concatenations for $A$ and $B$). See Figure \ref{BA_matrices} for an illustration.

Decomposing $A$, $B$, $D$, and $D'$ as compositions of $L$ and $R$ matrices, one can see that the composition of $A^{-1}$ and $B$ is conjugate to the composition of $D$ and $D'$:
\begin {align*} 
    A^{-1}B &= L_T^{-1} (D_{T'}D_T) L_T = R_T^{-1} (D_TD_{T'}) R_T, \quad \text{and} \\
    BA^{-1} &= R_{T'} (D_{T'}D_T) R_{T'}^{-1} = L_{T'} (D_T D_{T'}) L_{T'}^{-1}
\end {align*}

Since we have a combinatorial interpretation of the entries of $DD'$ by Lemma \ref{D-paths}, 
we will consider $DD'$ later on when we wish to examine the eigenvalues of $BA^{-1}$.

\subsection{Zig-zag paths}
Given some orientation of the graph $\G$ (either by (E), (NW), or (SW)), a \emph{zig-zag path} is a directed
path which alternately turns left, then right, then left, then right, etc. 
Let us also use the notation $H_A$ and $H_B$ for the concatenation of the networks used to define $A$ and $B$, as in Figure \ref{BA_matrices}.
These are both unfoldings of the graph $\G$.
The following result specifies how these paths correspond to the eigenvalues of the matrices $A$, $B$, and $BA^{-1}$.

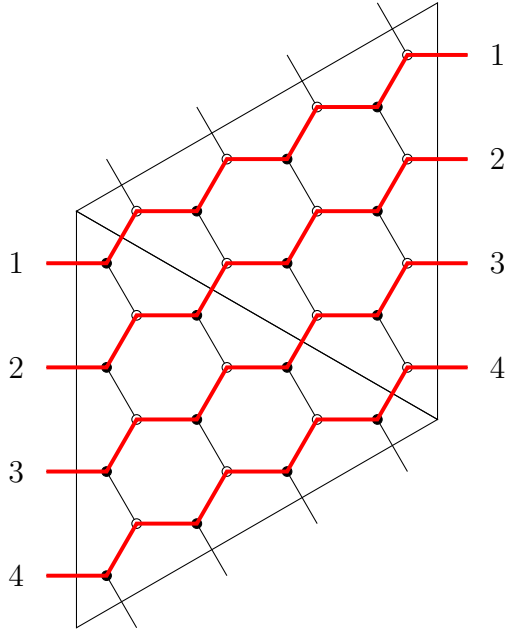
\begin {figure}[h]
\centering
\begin {tikzpicture}[scale=0.8]
    \draw (0:4) -- (120:4) -- (-120:4) -- cycle;
    \draw[shift={(60:4)}] (-60:4) -- (60:4) -- (180:4) -- cycle;

    \foreach \x/\y[evaluate={\X=int(\x+3); \Y=int(2-\y);}] in {0/0, 1/0, -1/0, -2/1, -1/1, 0/1, 1/1, 0/-1, 1/-1, 1/-2, -2/2, -1/2, 0/2, -2/3, -1/3, -2/4} {
        \draw ($({1.5*(\y+\x)},{0.5*sqrt(3)*(\y-\x)})$) -- ($({1.5*(\y+\x)},{0.5*sqrt(3)*(\y-\x)}) + (60:1)$);
        \draw ($({1.5*(\y+\x)},{0.5*sqrt(3)*(\y-\x)})$) -- ($({1.5*(\y+\x)},{0.5*sqrt(3)*(\y-\x)}) + (-60:1)$);
        \draw ($({1.5*(\y+\x)},{0.5*sqrt(3)*(\y-\x)})$) -- ($({1.5*(\y+\x)},{0.5*sqrt(3)*(\y-\x)}) + (180:1)$);
    }

    \begin {scope}[shift={(60:4)}]
        \foreach \y in {1,3,5,7} {
            \draw ($(-60:4) + (-0.5,{\y*sqrt(3)/2})$) -- ($(-60:4) + (0.5,{\y*sqrt(3)/2})$);
            \draw ($(60:4) + {\y*sqrt(3)/2}*(-150:1) + 0.5*(-60:1)$) -- ($(60:4) + {\y*sqrt(3)/2}*(-150:1) + 0.5*(120:1)$);
        }

        \draw ($(-60:4) + (1,{sqrt(3)/2})$) node {$4$};
        \draw ($(-60:4) + (1,{3*sqrt(3)/2})$) node {$3$};
        \draw ($(-60:4) + (1,{5*sqrt(3)/2})$) node {$2$};
        \draw ($(-60:4) + (1,{7*sqrt(3)/2})$) node {$1$};
    \end {scope}

    \foreach \x/\y in {0/0, 1/0, -1/0, -2/1, -1/1, 0/1, 1/1, 0/-1, 1/-1, 1/-2, -2/2, -1/2, 0/2, -2/3, -1/3, -2/4} {
        \draw[fill=black] ($({1.5*(\y+\x)},{0.5*sqrt(3)*(\y-\x)})$) circle (0.08);
    }
    \foreach \x/\y in {0/0, 1/0, -1/0, 0/-1, 1/-1, 1/-2, -2/1, -1/1, 0/1, 1/1, -2/2, -1/2, 0/2, -2/3, -1/3, -2/4} {
        \draw[fill=white] ($({1.5*(\y+\x)+0.5},{0.5*sqrt(3)*(\y-\x+1)})$) circle (0.08);
    }

    \draw[red, line width=1.5] ($(-120:4) + (-0.5,{7*sqrt(3)/2})$) --++ (0:1) --++ (60:1) --++ (0:1) --++(60:1) --++ (0:1) --++ (60:1) --++ (0:1) --++(60:1) --++(0:1);
    \draw[red, line width=1.5] ($(-120:4) + (-0.5,{5*sqrt(3)/2})$) --++ (0:1) --++ (60:1) --++ (0:1) --++(60:1) --++ (0:1) --++ (60:1) --++ (0:1) --++(60:1) --++(0:1);
    \draw[red, line width=1.5] ($(-120:4) + (-0.5,{3*sqrt(3)/2})$) --++ (0:1) --++ (60:1) --++ (0:1) --++(60:1) --++ (0:1) --++ (60:1) --++ (0:1) --++(60:1) --++(0:1);
    \draw[red, line width=1.5] ($(-120:4) + (-0.5,{1*sqrt(3)/2})$) --++ (0:1) --++ (60:1) --++ (0:1) --++(60:1) --++ (0:1) --++ (60:1) --++ (0:1) --++(60:1) --++(0:1);

    \draw ($(-120:4) + (-1,{sqrt(3)/2})$) node {$4$};
    \draw ($(-120:4) + (-1,{3*sqrt(3)/2})$) node {$3$};
    \draw ($(-120:4) + (-1,{5*sqrt(3)/2})$) node {$2$};
    \draw ($(-120:4) + (-1,{7*sqrt(3)/2})$) node {$1$};

\end {tikzpicture}
\caption{The zig-zag paths in the graph $H_A$.}
\label{A-zigzag}
\end {figure}

\begin {thm}\label{30casethm}
    Let $A$ and $B$ be defined as above. 
    \begin {enumerate}
        \item[(a)] The eigenvalues of $A$ are the weights of the north-east zig-zag paths in $\G$.
        \item[(b)] The eigenvalues of $B$ are the weights of the south-east zig-zag paths in $\G$.
        \item[(c)] The eigenvalues of $BA^{-1}$ are $(-1)^{n+1}$ times the weights of the south zig-zag paths in $\G$.
    \end {enumerate}
\end {thm}
\begin {proof}
    $(a)$ Recall that $A$ is a lower-triangular matrix, and so its eigenvalues are simply the diagonal entries. It is easy to see that with
    the (E) orientation, there is a unique path from $i$ on the left to $i$ on the right in $H_A$, given by the corresponding zig-zag path (see Figure \ref{A-zigzag}). 
    Thus the diagonal entry $A_{i,i}$ (and hence the eigenvalue) is the weight of this single path.

    $(b)$ Similarly, the diagonal entries of $B$ (which is upper-triangular) are the weights of unique zig-zag paths in $H_B$ from $i$ to $i$. 

    $(c)$ Recall that $BA^{-1}$ is conjugate to $DD'$, and so it has the same eigenvalues. By Lemma \ref{D-paths}, the entries of $DD'$
    are weighted sums of (SW) paths through $H_A$ from the top edge to the bottom edge. Essentially the same argument from $(a)$ and $(b)$
    works here, but with the following subtlety. 

    Although $D$ and $D'$ are themselves not upper or lower triangular, the product $DD'$ is lower triangular, as we will now explain.
    Because of the labels of the boundary vertices in the triangles $T$ and $T'$, one can see that $D$ is lower anti-triangular (i.e. $D_{ij} = 0$ for $i+j < n+1$)
    and $D'$ is upper anti-triangular (i.e. $D'_{ij} = 0$ for $i+j > n+1$). If $P$ is the permutation matrix for the permutation $k \mapsto n+1-k$,
    then $DP$ and $PD'$ are both lower triangular, and hence $(DP)(PD')$ is also lower-triangular. But since $P$ is an involution, $P^2 = \mathrm{Id}$, and so $(DP)(PD') = DD'$.

    Finally, since $DD'$ is lower-triangular, we have (as in parts $(a)$ and $(b)$) that its eigenvalues are the diagonal entries,
    which are (up to sign) the weights of the unique south zig-zag paths in the network between corresponding boundary vertices.
    However the signs in the product $DD'$ combine coherently to give total sign $(-1)^{n+1}$ to all entries in $DD'$.
\end {proof}

\section {Parameterization}\label{paramscn}

Recall that given multisets $\alpha = \{\alpha_1,\dots,\alpha_n\}$, $\beta = \{\beta_1,\dots,\beta_n\}$, and $\gamma = \{\gamma_1,\dots,\gamma_n\}$,
we consider the space $\Omega = \Omega_{0,3,n}(\alpha,\beta,\gamma)$ of pairs of matrices $(A,B)$ (considered up to conjugation; i.e. $(A,B) \cong (gAg^{-1},gBg^{-1})$)
such that $A$ has eigenvalues $\alpha_i$, $B$ has eigenvalues $\beta_i$, and $BA^{-1}$ has eigenvalues $\gamma_i$. The dimension of $\Omega$
is $(n-1)(n-2)$.

We can easily give an explicit description of $\Omega$ in the case $n=2$.

\begin {ex}
    For $n=2$, we have $\dim \Omega = 0$. The unique solution is represented by the edge weights in the network in Figure \ref{fig:n=2}.
    The matrices $A$ and $B$ are then given (up to conjugation) by     \[ 
        A = \begin{pmatrix} \alpha_1 & 0 \\ \alpha_1 - \frac{\beta_1}{\gamma_1} & \alpha_2 \end{pmatrix} \quad \quad \text{and} \quad \quad
        B = \begin{pmatrix} \beta_1 & \beta_2 - \alpha_2\gamma_1 \\ 0 & \beta_2 \end{pmatrix}
    \]
    Recall that $BA^{-1}$ is conjugate to $DD'$, which in this case is given by
    \[ DD' = \begin{pmatrix} \gamma_1 & 0 \\ \gamma_2 - \frac{\beta_2}{\alpha_2} & \gamma_2 \end{pmatrix} \]

\begin {figure}[h]
\centering
\begin {tikzpicture}[scale=1.2]
    \draw (0:2) -- (-120:2) -- (120:2);
    \draw[shift={(60:2)}] (-60:2) -- (60:2) -- (180:2);

    \foreach \x/\y in {0/0, 0/1, -1/1, -1/2} {
        \draw ($(-120:1) + ({1.5*(\y+\x)},{0.5*sqrt(3)*(\y-\x)})$) -- ($(-120:1) + ({1.5*(\y+\x)},{0.5*sqrt(3)*(\y-\x)}) + (60:1)$);
        \draw ($(-120:1) + ({1.5*(\y+\x)},{0.5*sqrt(3)*(\y-\x)})$) -- ($(-120:1) + ({1.5*(\y+\x)},{0.5*sqrt(3)*(\y-\x)}) + (-60:1)$);
        \draw ($(-120:1) + ({1.5*(\y+\x)},{0.5*sqrt(3)*(\y-\x)})$) -- ($(-120:1) + ({1.5*(\y+\x)},{0.5*sqrt(3)*(\y-\x)}) + (180:1)$);
    }

    \begin {scope}[shift={(60:2)}]
        \foreach \y in {1,3} {
            \draw ($(-60:2) + (-0.5,{\y*sqrt(3)/2})$) -- ($(-60:2) + (0.5,{\y*sqrt(3)/2})$);
            \draw ($(60:2) + {\y*sqrt(3)/2}*(-150:1) + 0.5*(-60:1)$) -- ($(60:2) + {\y*sqrt(3)/2}*(-150:1) + 0.5*(120:1)$);
        }

        \draw ($(-60:2) + (1,{1*sqrt(3)/2})$) node {$2$};
        \draw ($(-60:2) + (1,{3*sqrt(3)/2})$) node {$1$};
    \end {scope}

    \foreach \x/\y in {0/0, 0/1, -1/1, -1/2} {
        \draw[fill=black] ($(-120:1) + ({1.5*(\y+\x)},{0.5*sqrt(3)*(\y-\x)})$) circle (0.08);
    }
    \foreach \x/\y in {0/0, 0/1, -1/1, -1/2} {
        \draw[fill=white] ($({1.5*(\y+\x)},{0.5*sqrt(3)*(\y-\x)})$) circle (0.08);
    }

    \draw ($(-120:2) + (-1,{1*sqrt(3)/2})$) node {$2$};
    \draw ($(-120:2) + (-1,{3*sqrt(3)/2})$) node {$1$};

    \draw ($(0,{sqrt(3)}) + (-90:0.5)$) node {$\alpha_1$};
    \draw ($(0,{sqrt(3)}) + (100:0.8)$) node {$\beta_2$};
    \draw ($(0,0) + (90:0.5)$) node {$\beta_1$};
    \draw ($(0,0) + (-90:0.5)$) node {$-\alpha_2 \gamma_1$};
    \draw ($(1.5,{0.5*sqrt(3)}) + (-90:0.5)$) node {$-\gamma_1^{-1}$};

\end {tikzpicture}
\caption{The unique point in $\Omega(\alpha,\beta,\gamma)$ for $n=2$.}
\label{fig:n=2}
\end {figure}
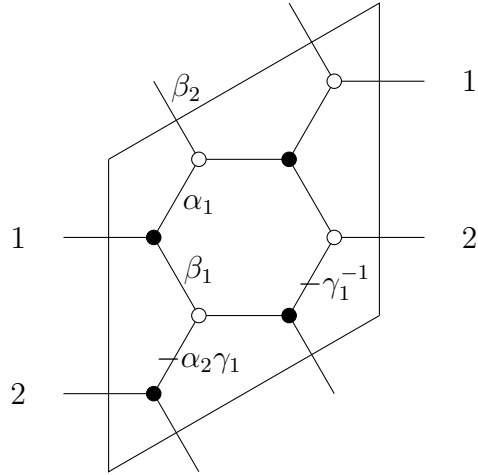
\end {ex}

\subsection{General $n$ case}

\begin {figure}[h]
\centering
\begin {tikzpicture}[scale=0.8, every node/.style={scale=0.8}]
    \foreach \x in {0,1,2,3} {
        \foreach \y in {0,1,2,3} {
        \draw ($({1.5*(\y+\x)},{0.5*sqrt(3)*(\y-\x)})$) -- ($({1.5*(\y+\x)},{0.5*sqrt(3)*(\y-\x)}) + (60:1)$);
        \draw ($({1.5*(\y+\x)},{0.5*sqrt(3)*(\y-\x)})$) -- ($({1.5*(\y+\x)},{0.5*sqrt(3)*(\y-\x)}) + (-60:1)$);
        \draw ($({1.5*(\y+\x)},{0.5*sqrt(3)*(\y-\x)})$) -- ($({1.5*(\y+\x)},{0.5*sqrt(3)*(\y-\x)}) + (180:1)$);
        \draw ($({1.5*(\y+\x)},{0.5*sqrt(3)*(\y-\x)}) + (-1,0)$) -- ($({1.5*(\y+\x)},{0.5*sqrt(3)*(\y-\x)}) + (-1,0) + (-120:1)$);
        \draw ($({1.5*(\y+\x)},{0.5*sqrt(3)*(\y-\x)}) + (-1,0)$) -- ($({1.5*(\y+\x)},{0.5*sqrt(3)*(\y-\x)}) + (-1,0) + (120:1)$);

        \draw ($({1.5*(\y+\x)},{0.5*sqrt(3)*(\y-\x)}) + (0:1)$) node {\scriptsize $X_{\x\y}$};
        }
    }
    \foreach \x in {0,1,2,3} {
        \foreach \y in {0,1,2,3} {
        \draw[fill=black] ($({1.5*(\y+\x)},{0.5*sqrt(3)*(\y-\x)})$) circle (0.08);
        \draw[fill=white] ($({1.5*(\y+\x)},{0.5*sqrt(3)*(\y-\x)}) + (-1,0)$) circle (0.08);
        }
    }

    \draw[blue,  line width=1.5] (-1.5,{-0.5*sqrt(3)}) --++ (60:1) --++ (0:1) --++ (60:1) --++ (0:1) --++ (60:1) --++ (0:1) --++ (60:1) --++ (0:1) --++ (60:1);
    \draw[red, line width=1.5] (150:{sqrt(3)}) --++ (-60:1) --++ (1,0) --++ (-60:1) --++ (1,0) --++ (-60:1) --++ (1,0) --++ (-60:1) --++ (1,0) --++(-60:1);
    \draw[purple, line width=1.5] (-1,0) --++ (1,0);
\end {tikzpicture}
\caption{Face weights in the honeycomb graph $H$, with the curves $\gamma_x,\gamma_y$ drawn in blue and red.}
\label{Xwts}
\end {figure}
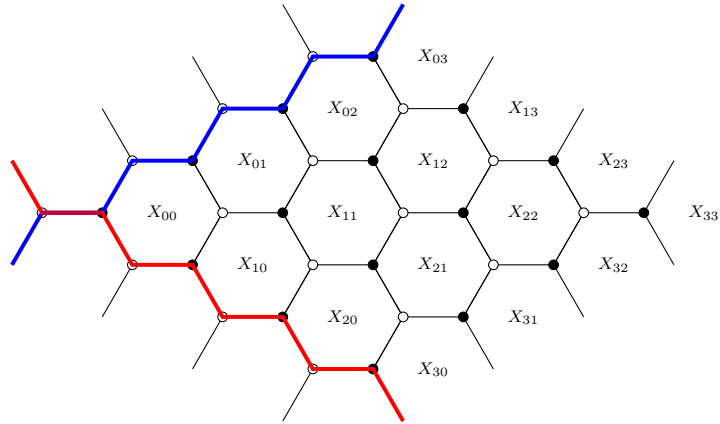

Given multisets $\alpha,\beta,\gamma$ of nonzero complex numbers, each of size $n$, choose an arbitrary order for each of them: 
$\alpha=(\alpha_1,\dots,\alpha_n)$ and likewise for $\beta$ and $\gamma$.

Let $\G_n$ be the torus honeycomb graph obtained by gluing $H_T$ and $H_{T'}$ (Figure \ref{A-zigzag})
and identifying opposite edges.
A $\C^*$-connection on $\G_n$ is determined by its monodromies along cycles in the cycle space of $\G_n$, which is of dimension $n^2+1$. 
The (counterclockwise) monodromies around faces $\{X_f\}_{f\in F}$ are called \emph{face weights}; 
there are $n^2$ of these quantities and there is one relation between them: their product is $1$.
There are additionally two extra cycles in the cycle space, generating the homology of the torus. Choose representative cycles $\gamma_x,\gamma_y$ and let
$X_x,X_y$ be their monodromies. For simplicity we can choose $\gamma_x,\gamma_y$ to be the (NE and SE, respectively) zigzag paths corresponding to $\alpha_1$ and $\beta_1$ (see Figure \ref{Xwts}). 

A $\C^*$-connection on $\G_n$ is determined up to gauge by specifying the values of $\XX = \{X_f: f\in F\}\cup\{X_x,X_y\}$.  
The eigenvalue equations determine the monodromies of the connection along $3n$ cycles, given by the $3n$ zig-zag paths. In \cite{KO}
it is shown that the space of $\C^*$-connections with given eigenvalues (that is, zig-zag path values) has dimension
$(n-1)(n-2)$, with face weights as global coordinates.
We give here an explicit and elementary proof. 

The eigenvalue equations are \emph{monomial} equations in the $\XX$ variables.
In particular the ratio of monodromies along two parallel zig-zag paths is the product of the face weights in the region between them:
for example $\alpha_1=X_x$ and $\alpha_2 = X_x\prod_{i=0}^{n-1} X_{0i}$. 
Likewise for each vertical zigzag path (oriented southwards) the monodromy is $X_y/X_x$ times the 
product of the $X_{i,j}^{-1}$ over all faces left of the path, using the fundamental domain of Figure \ref{Xwts}, and counting multiplicity.
For example for the zig-zag path consisting of the edges just right of the midline we have 
\be\label{midline}\gamma_1= \frac{X_y}{X_x}\prod_{i+j<n-1}X_{i,j}^{-1}.\ee

Fix $X_x=\alpha_1$ and $X_y=\beta_1$. 
We can now think of $\XX=(X_{i,j})_{i,j\in\{0,\dots,n-1\}}$ 
as an $n\times n$ matrix $M$ with given row products $\frac{\alpha_{i+1}}{\alpha_i}$, column products $\frac{\beta_{i+1}}{\beta_i}$ and ``large antidiagonal" products $\frac{\gamma_{i+1}}{\gamma_i}$, as well as one other relation, the midline equation (\ref{midline}) which is
$\prod_{i+j<n-1}X_{i,j}^{-1}=\frac{\gamma_1\alpha_1}{\beta_1}.$
We claim that the top left $(n-2)\times(n-1)$ block $\YY$ of entries of $M$ determines uniquely the remaining entries. This can be seen as follows.
Given the top $(n-2)\times(n-1)$ block of entries, let the first $n-1$ entries on the $(n-1)$th row be $p_1,\dots,p_{n-1}$
and those on the last row be $q_1,\dots,q_{n-1}$. (In the Figure \ref{Xwts}, $p_i=X_{2,i-1}$ and $q_i=X_{3,i-1}$.)
The entries on the last column (and first $n-2$ rows) are determined uniquely from $\YY$ by the row product condition.
The midline equation determines $p_1=X_{n-2,0}$ as a function of $\YY$ and $\alpha_1,\beta_1,\gamma_1$.
Once $p_1$ is determined, the column and antidiagonal product condition determine successively $q_1$, then $p_2$, then $q_2$, then $p_3$, and so on until $q_{n-1}$ is determined. Finally $p_n$ is determined by either the last row product or the last
antidiagonal product, using the fact that $\prod\frac{\alpha_i\gamma_i}{\beta_i}=1$. And $q_n$ is determined by any of the row, column or antidiagonal products. 

\subsection{Positivity}\label{positive}

In the case $n$ is odd, if the eigenvalues $\alpha,\beta,\gamma$ are all real and positive, then by Theorem \ref{30casethm},
there is a subset $\Omega^+_{0,3,n}\subset \Omega_{0,3,n}$ involving \emph{positive real} parameters $\XX$. 
This follows since each equation for the $\XX$ is monomial with coefficient $1$.
We thus have a positive Laurent parameterization
$$\Psi:(\R_+)^{(n-1)(n-2)}\to\Omega^+_{0,3,n}(\alpha,\beta,\gamma).
$$
Likewise if $n$ is even, and $\alpha,\beta$ consist in positive reals and $\gamma$ consists in negative reals,
then again we have a subset $\Omega^+_{0,3,n}\subset \Omega_{0,3,n}$ consisting of the image of the mapping with
positive real parameters.

\section {General $(g,k)$ }\label{genscn}\

In the previous sections, we discussed a parameterization of the space $\Omega_{0,3,n}(\lambda)$, when $(g,k) = (0,3)$, which was of dimension $(n-1)(n-2) = n^2-3n+2$.
We now discuss the general case of arbitrary $(g,k)$. The dimension of $\Omega_{g,k,n}(\lambda)$ in the general case will be $N=(2g-2+k)n^2-kn+2$,
and we will describe a parameterization analogous to the one given before, in terms of a graph embedded on a surface.

Let $\Sigma=\Sigma_{g,k}$ be a surface of genus $g$ with $k$ punctures $p_1,\dots,p_k$.
Fix $n$ and for each $i$ let $\{\lambda_{i,j}\}_{j=1}^n$ be $n$ eigenvalues associated to puncture $p_i$, satisfying $\prod_{i,j} \lambda_{i,j}=1$.

\begin{figure}[htbp]
\begin{center}\includegraphics[width=2.in]{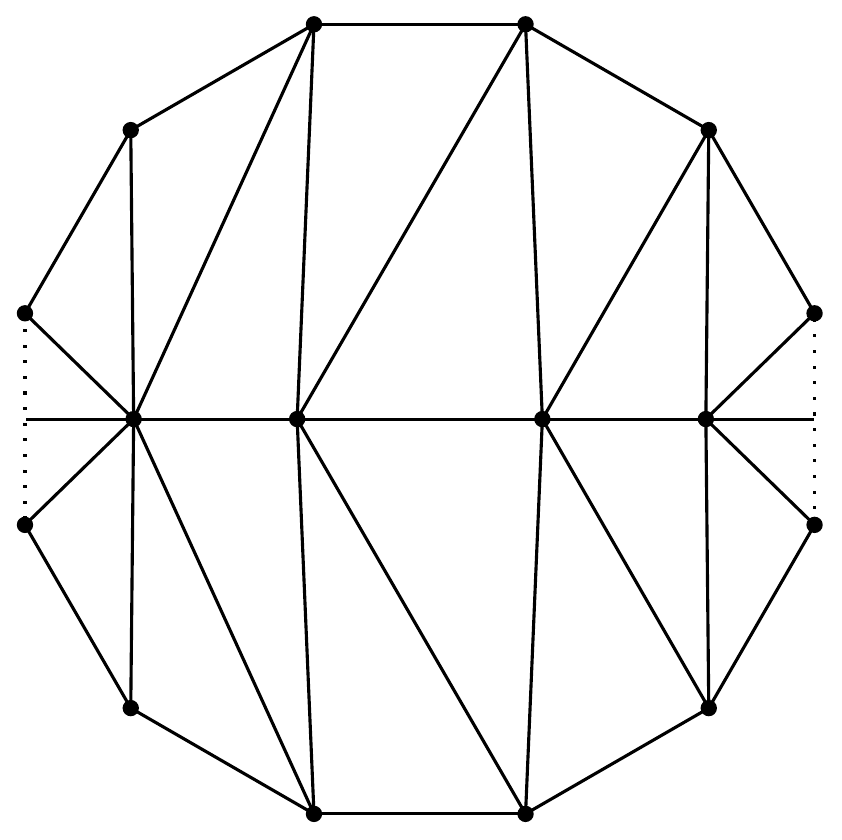}\end{center}
\caption{\label{triangulation}A genus-$g$ surface can be presented as shown as a regular $4g$-gon with opposite sides identified.
Take $k-1$ punctures along the horizontal symmetry axis and one puncture at the image of the vertices of the polygon.
(Shown is the case $g=3$ and $k=5$).
Take a triangulation containing the horizontal symmetry line and so that the triangulation of the upper half is the reflection of that in the lower half. Then all punctures will have even degree and so the dual graph of
the triangulation is bipartite.}
\end{figure}

The surface $\Sigma$ has a triangulation, with vertices at the punctures, with bipartite dual 
(see Figure \ref{triangulation}). 
There are $T$ triangles, where $T/2=2g-2+k$ is the negative of the Euler characteristic of $\Sigma$.

Let $\G_\Sigma$ be the dual graph of this triangulation.
The graph $\G_\Sigma$ has a natural ribbon structure (i.e. a cyclic ordering of the edges around each vertex) coming from its embedding in $\Sigma$.
There is a different ribbon structure on $\G_\Sigma$, which we denote $\G$, 
obtained by reversing the cyclic ordering at all black vertices.
It defines an embedding of $\G$ into another surface, denoted $S$. 
The surface $S$ has generally a different genus $g'$ and number of boundary components $k'$. 
However $S$ is triangulated with the same number of triangles $T$, so it has the same Euler characteristic as $\Sigma$, and thus $2-2g-k=2-2g'-k'$.

Zig-zag paths on $\G$ correspond to peripheral curves on $\Sigma$, and zig-zag paths
on $\G_\Sigma$ correspond to peripheral curves on $S$ 
(so $k$ is the number of zig-zag paths on $\G$).

\subsection{The graph $\G_n$}

Now fix $\G_\Sigma,\G,\Sigma,S$ as above. 
From $\G$ we make a new graph $\G_n$ (which we also draw embedded on the surface $S$) as follows.
In each triangle $t$ of $S$, corresponding to a white vertex of $\G$, we put a copy of
$H_t$ (see Figure \ref{Tn}). In each triangle $t'$ of $S$, corresponding to a black vertex of $\G$, we put a copy of
$H_{t'}$.

The zig-zag paths of $\G_n$ come in bands of $n$ parallel zig-zag paths 
associated to the same puncture of $\Sigma$. 

In the case $(g,k) = (0,3)$, discussed in earlier sections, the surface $\Sigma$ is a 3-holed sphere, and the conjugate surface $S$
is a once-punctured torus (i.e. $(g',k') = (1,1)$). This is why our honeycomb graph $H$ was drawn on a torus. 
The zig-zag paths in the three directions corresponding to $\alpha$, $\beta$, $\gamma$ correspond to the three holes of the sphere $\Sigma$.

\subsection{$\C^*$-connections}

For each puncture of $\Sigma$ there are $n$ parallel zig-zag paths of $\G_n$, and also $n$ prescribed eigenvalues. Pick an arbitrary bijection between these zig-zag paths and the eigenvalues.
We wish to define our $\C^*$-connection on $\G_n$ so that the monodromies of the zig-zag paths equal the
corresponding eigenvalues.

As previously, the $\C^*$-connection is defined by the set of face monodromies $\{X_f\}$ of $\G_n$
(including around punctures), with one condition that their product is $1$, and $2g'$ extra monodromies $\{X_{x_i},X_{y_i}\}_{i=1}^{g'}$ 
for a basis for the homology of the
closed surface $\bar S$. Thus 
the space of connections is parameterized by $\XX=\{X_f: f\in F\}\cup\{X_{x_i},X_{y_i}\}_{i=1}^{g'}$,
with the single relation above. 

The monodromy conditions for the zig-zag paths
give a set of $kn-1$ monomial equations for the $\XX$: there is one equation
for each eigenvalue/zig-zag path, but the last will be determined by the condition that the product of the eigenvalues is $1$. 

The cycle space of $\G_n$ is of dimension $F+2g'-1$ where $F$ is the number of its faces. We have $F=(n^2-1)T/2+k'$ where $T=4g+2k-4$ is the number of triangles in the original triangulation of $\Sigma$. This leads to the value
$$(n^2-1)T/2+k'+2g'-1=n^2T/2+1=n^2(2g-2+k)+1$$
for the dimension of the cycle space of $\G_n$. 

It remains to show that the $kn-1$ monomial equations are independent. 

Consider the conjugate surface $\widetilde{\G_n}$ of $\G_n$, obtained reversing the ribbon structure at each black vertex of $\G_n$. Zig-zag paths on $\G_n$ become peripheral curves on $\widetilde{\G_n}$,
which therefore are independent in the cycle space of $\G_n$, modulo their product being the
identity. This shows that the $kn-1$ monomial equations for the eigenvalues are independent. 

In conclusion the space of $\C^*$-connections on $\G_n$, after fixing the 
zig-zag monodromies, is $(\C^*)^N$ where
$$N=n^2(2g-2+k)-kn+2.$$

\section{Torus case and integrability}\label{intscn}

For certain $(g,k)$ and certain triangulations of $\Sigma_{g,k}$, the conjugate
surface $S$ will be a torus (with punctures). 
These are cases arising from a convex integer polygon $N\subset\Z^2$,
the \emph{Newton polygon}, as follows. 
Given a convex polygon $N$ with vertices in $\Z^2$, \cite{GK} constructed
a ``minimal'' bipartite graph $\G$ on the torus with Newton polygon $N$. Minimality is a technical assumption on the intersections of the zig-zag paths, see \cite{GK}, and by Newton polygon $N$ we mean that
the torus homology classes (in $\Z^2$) of zig-zag paths in $\G$, when arranged in order of increasing argument, form the primitive edges of $N$.
They showed that without loss of generality $\G$ can be chosen to have vertices of degree $2$ and $3$ only.

In this case the dual graph of $\G$ is not necessarily a triangulation, but has faces which are triangles and bigons. The graph $\G_n$ is constructed as before in each triangle but for a
black (resp. white) bigon face we simply put a sequence of $n$ parallel edges, each with a black
(resp. white) vertex in its center. 

The graph $\G_n$ is also
embedded and minimal on the torus $S$. 
In \cite{GK} it was proved (or one can repeat the argument of the previous section to show) that the zig-zag path monodromies of $\G_n$ on $S$ are free parameters on 
the space of connections (except for the one relation that their product is $1$).
Moreover they showed that these monodromies freely generate
the kernel of the associated Poisson structure, that is, define the space of Casimirs. 
The common level sets of the zig-zag path monodromies therefore has a symplectic structure. An explicit Hamiltonian
Liouville integrable system was constructed in \cite{GK}. 

It remains to be seen for which pairs $(g,k)$ there exists such a polygon $N$, that is, when is $S$ a punctured torus?
This is provided by a theorem of Scott:
 
\begin{prop}[\cite{Scott}]
For a convex polygon with integer vertices and having $g$ interior lattice points,
the number $k$ of integer boundary points satisfies $3\le k\le2g+6$, unless $g=1$ in which case $3\le k\le 9$, or $g=0$ in which case $3\le k$. 
 Moreover all such pairs $(g,k)$ are realized  by integer polygons.
\end{prop}

\begin{proof}
We provide a proof here of the second statement, which is what we need.
For the case $g=0$, take the triangle with vertices $(0,0),(k-2,0),(0,1)$.
For $(g,k)=(1,9)$ take the triangle with vertices $(0,0),(3,0),(0,3)$.
For $g>0$ and $(g,k)=(g,2g+6)$ we can take the convex polygon $N$ with vertices $(0,0),(0,2),(2g+2,0)$.
Now fixing $g$ we can decrease $k$ to any number $\ge 4$ by intersecting
$N$ with the upper half-space defined by the line through $(0,1)$ and $(i,0)$ for $1\le i\le 2g+2$. 
For the case
$k=3$ take $N$ with vertices $(0,1),(0,2),(2g+1,0)$.
\end{proof}

\bibliographystyle{alpha}
\bibliography{eig}

\end {document}